\documentclass[11pt]{article}

\usepackage{latexsym}
\usepackage{amsfonts}
\usepackage{enumerate}
\usepackage{multicol}
\usepackage{graphicx}
\usepackage{amssymb}
\usepackage{amsmath}
\usepackage{epic}

\begin{document}
\parindent=0.2in
\parskip 0in

\begin{flushright}

{\huge {\bf Report on the 49th Annual United States of America Mathematical Olympiad}}

\vspace*{.2in}

{\Large B\'ELA BAJNOK} \\

{\small Gettysburg College \\  Gettysburg, PA 17325 \\ bbajnok@gettysburg.edu}

\end{flushright}

\vspace*{.2in}

The American Mathematics Competitions program of the Mathematical Association of America consists of a series of examinations for middle school, high school, and college students, designed to build problem solving skills, foster a love of mathematics, and identify and nurture talented students across the United States.  The final round at the high-school level is the USA Mathematical Olympiad (USAMO).  This competition follows the style of the International Mathematics Olympiad: it consist of three problems each on two consecutive days, with an allowed time of four and a half hours both days.  

The 49th annual USA Mathematical Olympiad was given on Friday, June 19, 2020 and Saturday, June 20, 2020.  About 230 students were invited to take the exam; for the first time in history, the competition was administered online.  The names of winners and honorable mentions, as well as more information on the American Mathematics Competitions program, can be found on the site https://www.maa.org/math-competitions.  

The problems of the USAMO are chosen -- from a large collection of proposals submitted for this purpose -- by the USAMO Editorial Board, whose co-editors-in-chief are Evan Chen and Jennifer Iglesias, with associate editors John Berman, Zuming Feng, Sherry Gong, Alison Miller, Maria Monks Gillespie, and Alex Zhai.  This year's problems were created by Ankan Bhattacharya, Antonia Bluher, Zuming Feng, Carl Schildkraut, David Speyer, Richard Stong, and Alex Zhai.  

The solutions we present here are composed by the present author, and are based on the competition papers of Ankit Bisain (11th grade, Canyon Crest Academy, CA), Jeffrey Kwan (12th grade, Harker Upper School, CA), Rupert Li (12th grade, Jesuit High School, OR), Holden Mui (11th grade, Naperville North High School, IL), Yuru Niu (12th grade, Suncoast High School, FL), Ishika Shah (12th grade, Cupertino High School, CA), and Brandon Wang (12th grade, Saratoga High School, CA).  Each problem was worth 7 points; the nine-tuple $(n; a_7, a_6, a_5, a_4, a_3, a_2, a_1, a_0)$ states the number of students who submitted a paper for the relevant problem, followed by the numbers who scored $7, 6, \dots, 0$ points, respectively.

\vspace{.1in}

\noindent {\bf Problem 1} $(186; 143, 4, 0, 0, 0, 8, 16, 15)$ Let $ABC$ be a fixed acute triangle inscribed in a circle $\omega$ with center $O$.  A variable point $X$ is chosen on minor arc $AB$ of $\omega$, and segments $CX$ and $AB$ meet at $D$.  Denote by $O_1$ and $O_2$ the circumcenters of triangles $ADX$ and $BDX$, respectively.  Determine all points $X$ for which the area of triangle $OO_1O_2$ is minimized.  

\vspace{.1in}

\noindent {\em Solution.}  Let points $M$, $N$, $R$, and $S$ denote, in order, the midpoints of $\overline{CX}$, $\overline{DX}$, $\overline{AD}$, and $\overline{DB}$.  Since $C$ and $X$ are equidistant from $O$, $\overline{OM}$ is the perpendicular bisector of $\overline{CX}$. Similarly, $\overline{RO_1}$ is the perpendicular bisector of $\overline{AD}$, $\overline{SO_2}$ is the perpendicular bisector of $\overline{DB}$, and $\overline{O_1N}$ and $\overline{O_2N}$ are both perpendicular bisectors of $\overline{DX}$; in particular, $O_1$, $N$, and $O_2$ are collinear.  We then see that $RS = \frac{1}{2} AB$ and
$MN = MX - NX = \frac{1}{2} CX  - \frac{1}{2} DX = \frac{1}{2} CD.$

\begin{center}
	\includegraphics{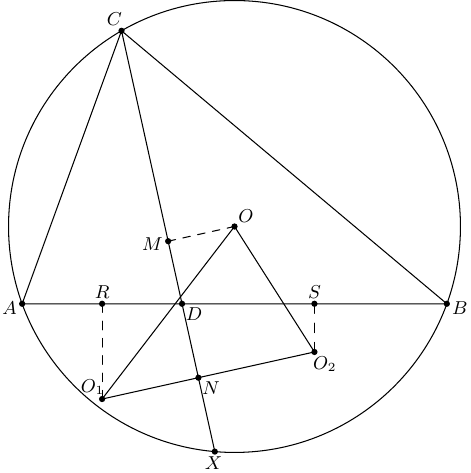}
\end{center}

Since in quadrilateral $RDNO_1$ the angles at $D$ and $O_1$ are supplementary, the angles $\angle ADC$ and $\angle RO_1O_2$ are equal; therefore, $RS = O_1O_2\sin \angle ADC$ and thus $$O_1O_2 = \frac{AB}{2 \sin \angle ADC}.$$  Furthermore, since lines $OM$ and $O_1O_2$ are parallel, the distance from $O$ to $\overline{O_1O_2}$ equals $MN$.  

We can now see that the area of triangle $OO_1O_2$ equals  
$$\frac{1}{2} \cdot O_1O_2 \cdot MN =  \frac{AB \cdot CD}{8 \sin \angle ADC}.$$  
Note that $CD$ is minimized when $\overline{CD} \perp \overline{AB}$, and $\sin \angle ADC$ is maximized when $\angle ADC = \pi/2$.  This is a happy coincidence, as the same point $X$ works for both of these conditions: the area of triangle $OO_1O_2$ is minimized when $X$ is the (unique) point with $\overline{CX} \perp \overline{AB}$.  
{\hfill $\Box$}

\vspace{.1in}

\noindent {\bf Problem 2}  $(156; 57, 35, 8, 2, 3, 9, 7, 35)$  An empty $2020 \times 2020 \times 2020$ cube is given, and a $2020 \times 2020$ grid of square unit cells is drawn on each of its six faces. A \emph{beam} is a $1 \times 1 \times 2020$ rectangular prism. Several beams are placed inside the cube subject to the following conditions:
\begin{itemize}
\item The two $1 \times 1$ faces of each beam coincide with unit cells lying on opposite faces of the cube. (Hence, there are $3 \cdot 2020^2$ possible positions for a beam.)
\item No two beams have intersecting interiors.
\item The interiors of each of the four $1 \times 2020$ faces of each beam touch either a face of the cube or the interior of the face of another beam.
\end{itemize}
What is the smallest positive number of beams that can be placed to satisfy these conditions?

\vspace{.1in}

\noindent {\em Solution.}  Let $n$ be a positive even integer, and consider an $n \times n \times n$ cube.  We claim that the smallest positive number of beams satisfying the three analogous conditions (where $2020$ is replaced by $n$) is $3n/2$ (and thus equals $3030$ for the case of $n=2020$).  

To facilitate our deliberation, we place the cube into the Cartesian coordinate system so that its edges are parallel to the coordinate axes and that one of its corners is at the origin and another at the point $(n,n,n)$.  We can then label the $n^3$ cells of the cube by their corners furthest from the origin.  For positive integers $a, b \leq n$, we then let $B_x(a,b)$ denote the {\em type $x$ beam} consisting of cells $(t,a,b)$ with $t = 1,2, \dots, n$; similarly, we let $B_y(a,b)$ denote the {\em type $y$ beam} consisting of cells $(a,t,b)$, and $B_z(a,b)$ denote the {\em type $z$ beam} consisting of cells $(a,b,t)$.

To see that there is a valid configuration with $3n/2$ beams, consider the collection consisting of beams
$$B_x(1,1), B_x(3,3), \dots, B_x(n-1,n-1),$$ $$B_y(1,n), B_y(3,n-2), \dots, B_y(n-1,2),$$ and $$B_z(2,2), B_z(4,4), B_z(6,6), \dots, B_z(n,n),$$ illustrated here for $n=10$.
 \begin{center}
	\includegraphics{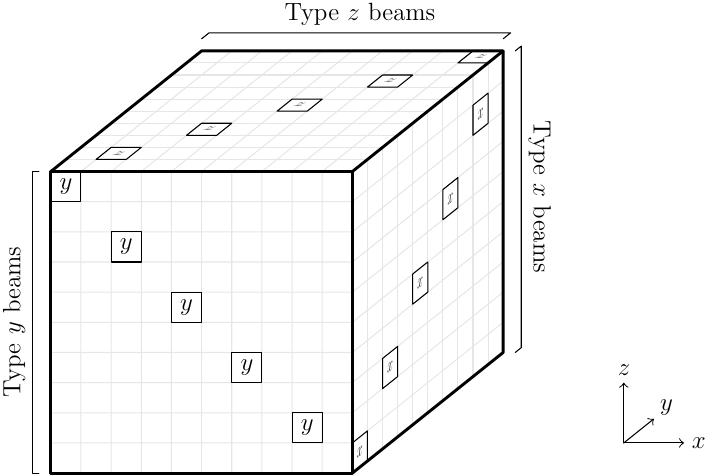}
\end{center}
Clearly, our collection consists of $3n/2$ beams satisfying the first requirement.

We can verify that our collection also satisfies the second condition, as follows.  Note that 
\begin{itemize}
  \item each type $x$ beam in the collection consists of cells whose second and third coordinates are odd, 
  \item each type $y$ beam consists of cells whose first coordinate is odd and third coordinate is even, and 
  \item each type $z$ beam consists of cells whose first and second coordinates are even.  
\end{itemize}  
Therefore, no two beams in the collection have intersecting interiors, and thus the second required condition holds.

Finally, to see that our third condition is satisfied, consider first a type $x$ beam $B$.  We see that the interior of each horizontal face (that is, each face perpendicular to the $z$ axis) of $B$ touches the interior of a type $y$ beam or the a face of the cube, and the interior of each vertical face (that is, each face perpendicular to the $y$ axis) of $B$ touches the interior of a type $z$ beam or the a face of the cube.  Since analogous observations apply to type $y$ beams and type $z$ beams, we conclude that our collection of beams satisfies all three required conditions.   

It remains to be shown that any valid construction must use at least $3n/2$ beams.  Let $N_x$, $N_y$, and $N_z$ denote the number of beams of type $x$, $y$, and $z$, respectively.  If any two of these three quantities are zero, then the third must equal $n^2$, since the collection must then contain all possible beams of that type.  As $n^2 > 3n/2$, we may assume that at least two of $N_x$, $N_y$, or $N_z$ are nonzero.

Next, we prove that $ N_x + N_y \geq n.$  By a {\em horizontal beam}, we mean a beam of type $x$ or type $y$; by the previous paragraph, we must have at least one horizontal beam.  According to the third condition, the interior of each horizontal face of each horizontal beam must touch the face of the cube or the interior of a face of another beam; clearly, this other beam would also need to be a horizontal beam.  Repeating this observation, we find that we must have at least $n$ horizontal beams, proving our claim.  

By a similar argument, we have $ N_y + N_z \geq n$ and $ N_z+N_x \geq n$, and  so 
$$ N_x + N_y + N_z = \frac12 (N_x+N_y) + \frac12 (N_y+N_z) + \frac12 (N_z+N_x)  \geq \frac32 n,$$
proving our claim.
{\hfill $\Box$}

\vspace{.1in}

\noindent {\bf Problem 3}  $(148; 30, 1, 1, 0, 5, 3, 11, 97)$  Let $p$ be an odd prime.  An integer $x$ is called a \emph{quadratic non-residue} if $p$ does not divide $x-t^2$ for any integer $t$.  

\noindent Denote by $A$ the set of all integers $a$ such that $1 \le a < p$, and both $a$ and $4-a$ are quadratic non-residues.  Calculate the remainder when the product of the elements of $A$ is divided by $p$. 

\vspace{.1in}

\noindent {\em Solution.}  We prove that the requested remainder is $2$.  As one can easily see that for $p=3$ we have  $A=\{2\}$, from here we will assume that $p \geq 5$, and work in the finite field $\mathbb{F}_p$, which here we identify with the set $\{0,1,2,\dots,p-1\}$ with operations carried out mod $p$. 

An element $n$ of $\mathbb{F}_p$ is a \emph{quadratic residue} when $n=t^2$ holds for some $t \in \mathbb{F}_p$ (it can be easily seen that there are at most two such $t$) and is a \emph{quadratic non-residue} otherwise (this corresponds to the definition as stated in the problem).  We will need the following facts: the product of two quadratic residues is a quadratic residue; the product of two quadratic non-residues is a quadratic residue; and the product of a quadratic residue and a quadratic non-residue is a quadratic non-residue.\footnote{For a proof of these claims and other properties of quadratic residues see, for example, Chapter VI in An introduction to the theory of numbers, by G.~H.~Hardy and E.~M.~Wright, Sixth edition, {\em Oxford University Press, Oxford}, 2008.}  

We let $A$ be the set of all elements $a \in \mathbb{F}_p$ such that both $a$ and $4-a$ are quadratic non-residues, and we let $B$ denote the set of all $b \in \mathbb{F}_p$ for which both $b$ and $4-b$ are quadratic residues. Note that $0 \in B$, $4 \in B$, and $2 \in A \cup B$.

Now consider an element $n \in A \cup B$.  Since $n(4-n)$ is then either the product of two quadratic residues or two quadratic non-residues, it must be a quadratic residue.  Furthermore, the element $4-n(4-n)$ is a quadratic residue as well, because it is equal to $(n-2)^2$.  This then implies that $n(4-n) \in B$.  

Conversely, we prove that for every $b \in B$, there is an element $n \in A \cup B$ such that $n(4-n)=b$.  Indeed, since $b \in B$, we have $4-b \in B$ as well, so there must be an element $n \in \mathbb{F}_p$ for which $(n-2)^2=4-b$, that is, $n(4-n)=b$.  But $b$ is a quadratic residue, so either both $n$ and $4-n$ are quadratic residues or they are both quadratic non-residues, so $n \in A \cup B$.  Moreover, we can see that, unless $b=4$, in which case $n=2$ is the only solution, the equation $n(4-n)=b$ has exactly two solutions, $n$ and $4-n$.

We thus see that there is a bijection between the unordered pairs $\{n,4-n\}$ (allowing for the pair $\{2,2\}$) of $A \cup B$ and the elements of $B$, given by the map $\{n,4-n\} \mapsto n(4-n)$.  Note that $\{0,4\}$ maps to $0$, and $\{2,2\}$ maps to $4$.  Consequently, when $2 \in A$, the product of the elements of $(A \cup B) \setminus \{0,2,4\}$ equals the product of the elements in $B \setminus \{0,4\}$.  This means that the product of the elements in $A \setminus \{2\}$ must be $1$, and thus the product of the elements in $A$ equals $2$.

The situation is similar when $2 \in B$, but this time the product of the elements of $(A \cup B) \setminus \{0,2,4\}$ equals twice the product of the elements in $B \setminus \{0,4\}$.  This implies that the product of the elements in $A$ again equals $2$.
{\hfill $\Box$} 

\vspace{.1in}

\noindent {\bf Problem 4}  $(179; 73, 1, 0, 0, 0, 0, 88, 17)$  Suppose that $(a_1, b_1), (a_2, b_2), \dots , (a_{100}, b_{100})$ are distinct ordered pairs
of nonnegative integers. Let $N$ denote the number of pairs of integers $(i, j)$ satisfying
$1 \leq i < j \leq 100$ and $|a_ib_j - a_jb_i| = 1$. Determine the largest possible value of $N$ over
all possible choices of the $100$ ordered pairs. 

\vspace{.1in} 

\noindent {\em First Solution.}  We claim that the answer is $N=2n-3$ for $n \geq 2$ ordered pairs ($197$ for $n=100$).  Let $P_1=(a_1,b_1), \dots, P_n=(a_n,b_n)$.  We say that points $P_i$ and $P_j$ are {\em enchanted} if $|a_ib_j - a_jb_i| = 1$; note that, by the Shoelace Formula, this is equivalent to triangle $OP_iP_j$ (where $O$ is the origin) having area $1/2$.  It is easy to see that the $n$ points $P_1=(0,1), P_2=(1,2), P_3=(1,3), \dots, P_n=(1,n)$ contain $2n-3$ enchanted pairs: $P_1$ is enchanted with the other $n-1$ points, and $P_i$ is enchanted with $P_{i+1}$ for each $i=2,3,\dots, n-1$.  

We will now use induction to prove that $N \leq 2n-3$ for every $n \geq 2$.  This being trivial for $n=2$, assume that our claim holds for each collection of $n-1$ points for some $n \geq 3$, and consider a collection $P_1=(a_1,b_1), \dots, P_n=(a_n,b_n)$.  Without loss of generality, we assume that $a_n+b_n \geq a_i+b_i$ for all $1 \leq i \leq n$.  By our inductive assumption, it suffices to show that $P_n$ is enchanted with at most two other points.  

If this were not the case, then we would have two points $P_i$ and $P_j$ that have the same distance from line $OP_n$ and are on the same side of that line.  But then $P_iP_j$ is parallel to $OP_n$, so $\overrightarrow{P_iP_j}=t \cdot \overrightarrow{OP_n} = \langle ta_n, tb_n \rangle$ for some scalar $t$, and we may assume that $t>0$.  Since $P_i$ and $P_j$ have integer coordinates, $ta_n$ and $tb_n$ are integers, so $|b_i ta_n - a_i tb_n|=t$ is an integer as well.  Therefore, $t \geq 1$, and since $P_i \neq O$, we arrive at $a_j+b_j=(a_i+ta_n)+(b_i+tb_n) > a_n+b_n$, contradicting our choice of $P_n$.  
{\hfill $\Box$}  

\vspace{.1in} 

\noindent {\em Second Solution.}  First, we recall that the {\em Farey sequence} of order $m$ is the increasing sequence of fractions $a/b$ of relatively prime integers $a$ and $b$ with $0 \leq a \leq b \leq m$.  The property of Farey sequences that we need here is that $a/b$ and $a'/b'$ are consecutive terms in some Farey sequence if, and only if, $|a'b-ab'|=1$.\footnote{For this and other interesting features of Farey sequences see, for example, Chapter III in An introduction to the theory of numbers, by G.~H.~Hardy and E.~M.~Wright, Sixth edition, {\em Oxford University Press, Oxford}, 2008.}     

Suppose now that $(a_1, b_1), (a_2, b_2), \dots , (a_n, b_n)$ are distinct ordered pairs
of nonnegative integers, so that there are $N$ pairs of integers $(i, j)$ satisfying
$1 \leq i < j \leq n$ and $|a_ib_j - a_jb_i| = 1$. We shall use induction to prove that $N \leq 2n-3$ for every $n \geq 2$.  This obviously holds for $n=2$ and $n=3$, so let $n \geq 4$.  Without loss of generality, we arrange our points so that $\max (a_n, b_n) \geq \max (a_i,b_i)$ for each $1 \leq i \leq n$; furthermore, we may assume that $a_n \leq b_n$, since if this were not the case, we could instead consider the mirror image of our $n$ points with respect to the line $y=x$.  Note that our assumptions imply that $b_n \geq 2$, since we can only have $b_n=1$ for $n \leq 3$.  Our goal is to show that there are at most two indices $1 \leq i \leq n-1$ for which $|a_ib_n - a_nb_i| = 1$; our claim will then follow by induction.  

Observe that if $|a_ib_n - a_nb_i| = 1$, then $a_i \leq b_i$.  Indeed, if this were not the case, then we would have $a_i \geq b_i+1$, so     
$$1=|a_ib_n - a_nb_i| = a_ib_n - a_nb_i \geq (b_i+1)b_n-a_nb_i = b_i(b_n-a_n) +b_n \geq b_n \geq 2,$$ a contradiction.  Note also that $|a_ib_n - a_nb_i| = 1$ implies that $b_i \neq 0$ and  $\gcd (a_i,b_i)=\gcd (a_n,b_n)=1$.  Therefore, $a_i/b_i$ and $a_n/b_n$ are consecutive terms in some Farey sequence; we need to show that this can happen for at most two indices $i$.  But this is clearly the case: $a_n/b_n$ appears only in Farey sequences of order $b_n$ or more; it can have at most two neighbors in the Farey sequence of order $b_n$; and if $a_i/b_i$ and $a_n/b_n$ are not consecutive in the Farey sequence of order $b_n$, then they are also not consecutive in one with order more than $b_n$.

To show that $N=2n-3$ is achievable, we may start with any two fractions $a/b$ and $a'/b'$ that are consecutive in some Farey sequence; inserting $(a+a')/(b+b')$ between them results in three consecutive fractions, since
$|b(a+a')-a(b+b')|=1$ and $|b'(a+a')-a'(b+b')|=1$.  Repeating this process then yields $n$ points with $N=2n-3$ pairs with the desired property.  (For example, starting with $0/1$ and $1/2$, and inserting the just-described fraction next to $0/1$ each time, yields the sequence $0/1, 1/n, 1/(n-1), \dots, 1/2$, corresponding to the example of the first solution.)
{\hfill $\Box$} 

\vspace{.1in}

\noindent {\bf Problem 5} $(140; 21, 9, 1, 9, 17, 17, 49, 17)$  A finite set $S$ of points in the coordinate plane is called {\em overdetermined} if $|S| \geq 2$ and there exists a nonzero polynomial $P(t)$, with real coefficients and of degree at most $|S|-2$, satisfying $P(x)=y$ for every point $(x,y) \in S$.  For each integer $n \geq 2$, find the largest integer $k$ (in terms of $n$) such that there exists a set of $n$ distinct points that is {\em not} overdetermined, but has $k$ overdetermined subsets.  

\vspace{.1in}

\noindent {\em Solution.}  Recall that for every nonempty finite set $S$ in the coordinate plane consisting of points with distinct $x$ coordinates, there exists a unique polynomial $f_S$, called the {\em interpolating polynomial of $S$}, which has real coefficients and  degree at most $|S|-1$, and which satisfies $f(x)=y$ for every point $(x,y) \in S$.  Thus we can say that $S$ is overdetermined if, and only if, it has at least two elements and its interpolating polynomial has degree at most $|S|-2$.  (Note that sets containing points that share their $x$ coordinates do not have interpolating polynomials.)

Let $n \geq 2$, and consider the set $$A=\{(1,2)\} \cup \{(2,1), (3,1), \dots, (n,1)\}.$$  Then $$f_A(x)=1+\frac{1}{(n-1)!} (2-x)(3-x) \cdot \cdots \cdot (n-x),$$ so $A$ is not overdetermined; however, each of its subsets not containing $(1,2)$ and having size at least $2$ has interpolating polynomial of degree $0$ and is thus overdetermined.  So we found a non-overdetermined set of size $n$ with at least
\begin{eqnarray} \label{USAMO 5 identity}
{n-1 \choose 2} + {n-1 \choose 3} + \cdots + {n-1 \choose n-1} & = & 2^{n-1}-n
\end{eqnarray} overdetermined subsets; we will prove that we cannot do better.

Let $A$ be a set of $n \geq 2$ distinct points in the coordinate plane; we will prove that if it has more than $2^{n-1}-n$ overdetermined subsets, then it is overdetermined.  Let $N_m$ denote the number of overdetermined $m$-subsets of $A$; by (\ref{USAMO 5 identity}) above, it suffices to show that if $N_m > {n-1 \choose m}$ for some $2 \leq m \leq n-1$, then $A$ is overdetermined.  We will, in fact, prove the stronger statement that if $N_m > {n-1 \choose m}$ for some $2 \leq m \leq n-1$, then not just $A$ but all its subsets of size at least $m$ are overdetermined.

Keeping $m$ as fixed, we use induction on $n$.  Suppose first that $n=m+1 \geq 3$ and that $N_m >1$.  Let $P_1=(x_1,y_1)$ and $P_2=(x_2,y_2)$ be distinct points in $A$ for which $A \setminus P_1$ and $A \setminus P_2$ are overdetermined.  Then $f_{A \setminus P_1}$  and $f_{A \setminus P_2}$ are polynomials of degree at most $|A|-3$; since they agree on $|A|-2$ values of $x$, they must equal for all values of $x$, and thus $f_{A \setminus P_1}(x_1)=f_{A \setminus P_2}(x_1)=y_1$ and $f_{A \setminus P_2}(x_2)=f_{A \setminus P_1}(x_2)=y_2$.  This means that $A$ has an interpolating polynomial and it is of degree at most $|A|-3$, which then implies that $A$ and all its subsets of size $|A|-1$ are overdetermined, as claimed.

Suppose now that our claim holds for all sets of size $n-1$, and consider a set $A$ of size $n$ with $N_m > {n-1 \choose m}$ for some $2 \leq m \leq n-1$.  Since each $m$-subset of $A$ is in $n-m$ subsets of size $n-1$ and since
$$\frac{{n-1 \choose m}}{{n-2 \choose m}}=\frac{n-1}{n-m-1} > \frac{n}{n-m},$$ at least one of the $n$ subsets of $A$ of size $n-1$ has more than ${n-2 \choose m}$ overdetermined subsets of size $m$; let one such subset be $A \setminus \{P\}$.  By our inductive assumption, all subsets of $A \setminus \{P\}$ of size at least $m$ are overdetermined.  Since $N_m > {n-1 \choose m}$, there is an overdetermined set of $A$ of size $m$ that is not a subset of $A \setminus \{P\}$; let $S$ be one of these sets.  Let $Q \not \in S$, and consider $A \setminus \{Q\}$.

Note that the ${n-2 \choose m}$ $m$-subsets of $A \setminus \{Q, P\}$ are all overdetermined, and so is $S$.  Therefore, by our inductive hypothesis again, all subsets of $A \setminus \{Q\}$ of size at least $m$ are overdetermined.  By now we have at least two overdetermined subsets of $A$ of size $n-1$, so $A$ is overdetermined as well.

It remains to be shown that all subsets of $A$ containing both $P$ and $Q$ and having size between $m$ and $n-1$, inclusive, are overdetermined.  For this purpose, we let $R$ be an arbitrary element of $A \setminus \{P, Q\}$, and show that all subsets of $A \setminus \{R\}$ of size at least $m$ are overdetermined.  Note that all $m$-subsets of $A \setminus \{P, R\}$ and all $m$-subsets of $A \setminus \{Q, R\}$ are overdetermined, giving a total of more than ${n-2 \choose m}$ overdetermined $m$-subsets of $A \setminus \{R\}$, so our claim follows from our inductive hypothesis.  Our proof is now complete.
{\hfill $\Box$} 

\vspace{.1in}

\noindent {\bf Problem 6}  $(88; 4, 2, 1, 0, 0, 0, 1, 80)$   Let $n \geq 2$ be an integer.  Let $x_1 \geq x_2 \geq \cdots \geq x_n$ and $y_1 \geq y_2 \geq \cdots \geq y_n$ be $2n$ real numbers such that
$$0=x_1 + x_2 + \cdots + x_n = y_1 + y_2 + \cdots + y_n$$ and $$1=x_1^2 + x_2^2 + \cdots + x_n^2 = y_1^2 + y_2^2 + \cdots + y_n^2.$$
Prove that
$$\sum_{i=1}^n (x_iy_i-x_iy_{n+1-i}) \geq \frac{2}{\sqrt{n-1}}.$$  

\vspace{.1in}

\noindent {\em Solution.}  We will use the Rearrangement Inequality, which can be stated as follows.  Suppose that $x_1 \geq x_2 \geq \cdots \geq x_n$ and $y_1 \geq y_2 \geq \cdots \geq y_n$ are real numbers.  Let $S_n$ be the set of permutations of $\{1,2,\dots,n\}$, and for each $\sigma \in S_n$, set $f(\sigma) = \sum_{i=1}^n x_i y_{\sigma(i)}.$  With these notations, $f(\sigma)$ is maximized when $\sigma$ is the identity permutation (that is, $\sigma (i)=i$ for each $1 \leq i \leq n$) and minimized when $\sigma$ is the reverse identity permutation (that is, $\sigma (i)=n+1-i$ for each $1 \leq i \leq n$).\footnote{See Chapter 10 in Inequalities, by G.~H.~Hardy, J.~E.~Littlewood, and G.~P\'olya, Second Edition, {\em Cambridge University Press, London}, 
1952.}   

We will also need the following.

{\bf Lemma.}  Suppose that $a_1, \dots, a_k$ are real numbers, with $a_1 \geq a_i \geq a_k$ for all $1 \leq i \leq k$.  Suppose that $0 = \sum_{i=1}^k a_i$ and $A = \sum_{i=1}^k a_i^2$.  Then 
$$a_1-a_k \geq 2 \sqrt{A/k}.$$

{\em Proof of Lemma}.  Let $a=(a_1+a_k)/2$.  Then
$$\sum_{i=1}^k (a_i-a)^2= \sum_{i=1}^k a_i^2 - 2 a \sum_{i=1}^k a_i + \sum_{i=1}^k a^2 \geq A.$$ Note that $(a_1-a)^2=(a_k-a)^2 \geq (a_i-a)^2$ for all $1 \leq i \leq k$, so $a_1-a\geq \sqrt{A/k}$ and $a-a_k \geq \sqrt{A/k}$; adding these two inequalities proves our claim.

We will now compute $\sum_{\sigma \in S_n} f(\sigma)$ and $\sum_{\sigma \in S_n} f^2(\sigma)$ under the conditions of this problem.  The first of these is easy: since $\sum_{i=1}^n x_i=0$, we have
$$\sum_{\sigma \in S_n} f(\sigma) = \sum_{\sigma \in S_n}\sum_{i=1}^n x_i y_{\sigma(i)}  = \sum_{i=1}^n x_i  \sum_{\sigma \in S_n} y_{\sigma(i)} = 0.$$

In order to compute $\sum_{\sigma \in S_n} f^2(\sigma)$, we first rewrite it as 
$$\sum_{\sigma \in S_n} \left(\sum_{i=1}^n x_i y_{\sigma(i)} \right)^2 = 
\sum_{\sigma \in S_n} \left( \sum_{i=1}^n x_i^2 y_{\sigma(i)}^2  + \sum_{i \neq j } x_i x_j y_{\sigma(i)} y_{\sigma(j)}\right) = \sum_{i=1}^n x_i^2 \sum_{\sigma \in S_n} y_{\sigma(i)}^2 + \sum_{i \neq j } x_i x_j \sum_{\sigma \in S_n} y_{\sigma(i)} y_{\sigma(j)}.$$ 
Here $\sum_{i=1}^n x_i^2=1$, and 
$$\sum_{i \neq j } x_i x_j =  \left(\sum_{i=1}^n x_i \right)^2 - \sum_{i=1}^n x_i^2 = 0^2-1=-1.$$
Next, note that for each $1 \leq i, i' \leq n$ there are exactly $(n-1)!$ permutations $\sigma \in S_n$ for which $\sigma (i)=i'$, so 
$$\sum_{\sigma \in S_n} y_{\sigma(i)}^2 = (n-1)! \sum_{i'=1}^n y_{i'}^2=(n-1)! \cdot 1.$$
Similarly, observe that for each $1 \leq i, i',j,j' \leq n$ with $i \neq j$ and $i' \neq j'$ there are exactly $(n-2)!$ permutations $\sigma \in S_n$ for which $\sigma (i)=i'$ and $\sigma (j)=j'$, so
$$\sum_{\sigma \in S_n} y_{\sigma(i)} y_{\sigma(j)}=(n-2)! \sum_{i' \neq j' } y_{i'} y_{j'} = (n-2)! \cdot (-1).$$
In summary, we find that 
$$\sum_{\sigma \in S_n} f^2(\sigma) = 1 \cdot (n-1)! \cdot 1 + (-1) \cdot (n-2)! \cdot (-1) = n(n-2)!.$$

Now let $\sigma_1, \sigma_2, \dots, \sigma_{n!}$ be the elements of $S_n$ in some order so that $\sigma_1$ is the identity permutation and $\sigma_{n!}$ is the reverse identity permutation.  Set $k=n!$ and $a_i=f(\sigma_i)$.  Then, by the Rearrangement Inequality, we have $f(\sigma_{n!}) \leq f(\sigma_i) \leq f(\sigma_1)$, so our Lemma yields
$$f(\sigma_1)-f(\sigma_{n!}) \geq 2 \sqrt{n(n-2)!/n!}$$ or  
$$\sum_{i=1}^n (x_iy_i-x_iy_{n+1-i}) \geq \frac{2}{\sqrt{n-1}},$$ as claimed.
{\hfill $\Box$}

\vspace{.1in}

\noindent {\scriptsize {\bf Acknowledgments.} The author wishes to express his immense gratitude to everyone who contributed to the success of the competition: the students and their teachers, coaches, and proctors; the problem authors; the USAMO Editorial Board; the graders; the Art of Problem Solving; and the AMC Headquarters of the MAA.  I am also grateful to Evan Chen for proofreading this article and for providing the two figures.} 

\vspace*{.3in}

\noindent {\bf B\'ELA BAJNOK} (MR Author ID: 314851) is a Professor of Mathematics at Gettysburg College and the Director of the American Mathematics Competitions program of the MAA.

\end{document}